\documentclass[11pt]{article}

\usepackage{amssymb}
\usepackage{amsmath}
\usepackage{amsthm}
\usepackage{amsfonts}

\title{\bf The relation of semiadjacency\\
on $\cap$-semigroups of transformations}
\author{Wieslaw A. Dudek and Valentin S. Trokhimenko\\[4pt]
{\small Communicated by M. V. Volkov}}
\date{}

\begin{document}
\sloppy \maketitle
\newtheorem{Theorem}{Theorem}
\newtheorem{Proposition}{Proposition}
\newtheorem{Definition}{Definition}
\newtheorem{Coll}{Corollary}

\newtheorem{theorem}{Theorem}
\newtheorem{proposition}{Proposition}
\newtheorem{lemma}{Lemma}
\newtheorem{definition}{Definition}
\newtheorem{corollary}{Corollary}
\newcommand{\pr}{\mbox{pr}_1}
\newcommand{\prr}{\mbox{pr}_2}
\newcommand{\spr}{\mbox{\scriptsize pr}_1}
\newcommand{\sprr}{\mbox{\scriptsize pr}_2}
\newcommand{\hi}{\chi_{\Phi}}
\newcommand{\shi}{\raisebox{-3pt}{$\begin{smallmatrix}{\small\chi}\\{\scriptscriptstyle
1}\end{smallmatrix}_{\hspace{-1mm}\scriptscriptstyle\Phi}$}}
\renewcommand{\leq}{\leqslant}

\begin{abstract}\noindent
We consider two relations on a $\cap$-semigroup of partial
functions on a given set: the inclusion of domains and
semiadjacency (i.e., the inclusion of the image of the first
function in the domain of the second). These are characterized
from an abstract point of view via a system of elementary axioms,
i.e., conditions expressed in the language of pure predicate
calculus with equality.

\medskip

\end{abstract}

\medskip

{\bf 1.} Let $\mathcal{F}(A)$ be the set of all transformations
(i.e., the partial maps) of a non-empty set $A$. The composition
of such maps is defined as $(g\circ f)(a)=g(f(a))$, where for every
$a\in A$ the left and right hand sides are defined, or undefined,
simultaneously (cf. \cite{Clif}). If a set $\Phi\subset\mathcal{F}(A)$
is closed under the composition, then the algebra $(\Phi,\circ)$
is called a \emph{semigroup of transformations} (cf. \cite{Sch3})
since the composition is an associative operation. If $\Phi$ is
also closed under the set-theoretic intersection $\cap$ of
transformations treated as subsets of $A\times A$, then the
algebra $(\Phi,\circ,\cap)$ is called a \emph{$\cap $-semigroup
of transformations}. The first abstract characterization of such
semigroups was found by Garvatski\u{\i} \cite{Garv}.
This characterization is based on the method of determining pairs
proposed by Schein \cite{Sch1}.

If $(G,\cdot)$ is a semigroup, then, as usual, $G^1$ denotes $G$
whenever $(G,\cdot)$ possesses an identity element and
$G^1=G\cup\{e\}$ otherwise; in the latter case the multiplication
$\cdot$ is extended to $G^1$ such that $e$ becomes an identity
element.

The main results obtained by Garvatski\u{\i} are presented
in the following two theorems.

\begin{theorem}\label{T-dt1} An algebra $(G,\cdot,\wedge)$ is isomorphic to some
$\cap$-semigroup of transformations if and only if $(G,\cdot)$ is
a semigroup, $(G,\wedge)$  is a semilattice and
\begin{eqnarray}
&&\label{dt-1} x(y\wedge z)=xy\wedge xz,\\
&&\label{dt-2} (x\wedge y\wedge z) u\wedge
(y\wedge z) v=(x\wedge y) u\wedge (y \wedge z)v
\end{eqnarray}
for all $x,y,z\in G$ and $u,v\in G^1$.
\end{theorem}

\begin{theorem}\label{T-dt2}
An algebra $(G,\cdot,\wedge)$ is isomorphic to some
$\cap$-semigroup of invertible $($i.e., one-to-one$)$
transformations if and only if $(G,\cdot)$ is a semigroup,
$(G,\wedge)$ is a semilattice, the operation $\cdot$
distributes over the operation $\wedge$ on the left and on
the right and the equality
\begin{equation}
\label{dt-4} xv\wedge uv\wedge uy\wedge xy=xv\wedge uv\wedge uy,
\end{equation}
is satisfied for all $x,y,u,v\in G^1$ such that
$xv,uv,uy,xy\in G$.
\end{theorem}

\medskip

{\bf 2.} Let $(G,\cdot)$ be a semigroup, $\rho $ a binary relation
on $G$, i.e., $\rho\subset G\times G$. Following \cite{Clif, Sch3}
we say that this relation is:
\begin{itemize}
\item \emph{stable} if   $(x,y),\, (u,v)\in\rho$ implies $(xu,yv)\in\rho$   for all $x,y,u,v\in G$;
\item \emph{left regular} if   $(u,v)\in\rho$ implies $(xu,xv)\in\rho$   for any $x,u,v\in G$;
\item \emph{right regular} if   $(x,y)\in\rho$ implies $(xu,yu)\in\rho$    for all $x,y,u\in G$;
\item \emph{left ideal} if $(x,y)\in\rho$ implies $(ux,y)\in\rho$   for any $x,y,u\in G$;
\item \emph{right negative} if  $(x,yu)\in\rho$ implies $(x,y)\in\rho$ for all $x,y,u\in G$.
\end{itemize}

A reflexive and transitive binary relation is called a \emph{quasi-order}.
It is clear that a quasi-order is stable if and only if it is left and right
regular. Also it is not difficult to see that a quasi-order $\rho$ on a semigroup
$(G,\cdot)$ is right negative if $(xy,x)\in\rho$ holds for all $x,y\in G$.

\medskip

{\bf 3.} Let $(\Phi,\circ)$ be a semigroup of transformations. We
will consider three relations $\zeta_{\Phi}$, $\chi_{\Phi}$ and
$\delta_{\Phi}$ on the set $\Phi$, defined as follows:
\begin{eqnarray}
  & &\zeta_{\Phi} =\{(f,g)\,|\,f\subset g\},\nonumber\\
  & &\chi_{\Phi}=\{(f, g)\,|\,\pr f\subset\pr g\},\nonumber\\
  & &\delta_{\Phi} =\{(f,g)\,|\,\prr f\subset\pr g\},\nonumber
\end{eqnarray}
where $\pr f$ and $\prr f$ denote the domain and the image of $f$,
respectively. The first relation is called the \emph{fundamental
order}, the second is the \emph{projection quasi-order}, the
third is the \emph{semiadjacency relation} on $\Phi$.

\begin{proposition}\label{P-trp1} For any $f,g\in\Phi$ we have
\begin{equation}\label{f-tr3}
(f,g)\in\delta_{\Phi}\longleftrightarrow (f,g\circ f)\in\hi.
\end{equation}
\end{proposition}
\begin{proof} Let $(f,g)\in\delta_{\Phi}$, i.e., $\prr f\subset\pr g$. Then
$f^{-1}(\prr f)\subset f^{-1}(\pr g)$, whence, in view of
$f^{-1}(\prr f)=\pr f$ and $f^{-1}(\pr g)=\pr (g\circ f)$, we
obtain $\pr f\subset\pr (g\circ f)$. Thus, $(f,g\circ f)\in\hi$.

Conversely, if $\pr f\subset\pr (g\circ f)$ and $b\in\prr f$, then
$(a,b)\in f$ for some $a\in A$. But $a\in\pr f$, so $a\in\pr
(g\circ f)$. Consequently, $(a,d)\in f$ and $(d,c)\in g$ for some
$d\in A$. Since $f$ is a transformation, from $(a,b)\in f$ and
$(a,d)\in f$ it follows $b=d$. Thus, $(b,c)\in g$, which implies
$b\in\pr g$. So, we have shown that $\prr f\subset\pr g$, that is,
$(f,g)\in\delta_{\Phi}$.
\end{proof}

\begin{proposition}\label{P-trp3} In a semigroup $(\Phi,\circ)$ the relation
$\delta_{\Phi}$ is left ideal.
\end{proposition}
\begin{proof} Let $f,g,u\in\Phi$ and $(f,g)\in\delta_{\Phi}$. Then
$\prr f\subset\pr g$. Since $\prr(f\circ u)\subset\prr f$, we have
$\prr (f\circ u)\subset\pr g$. Thus $(f\circ
u,g)\in\delta_{\Phi}$. So, $\delta_{\Phi}$ is left ideal.
\end{proof}

According to \cite{Sch3} semigroups of the form
$(\Phi,\circ,\zeta_{\Phi})$ are called \emph{ordered semigroups
of transformations}, semigroups of the form
$(\Phi,\circ,\zeta_{\Phi},\chi_{\Phi})$ are called
\emph{projection ordered semigroups of transformations}.
Abstract characterizations of these semigroups are given in
\cite{Sch3}. Semigroups of the form
$(\Phi,\circ,\cap,\chi_{\Phi})$ are called \emph{projection
$\cap$-semigroups of transformations}. These semigroups are a
special case of \emph{projection quasi-ordered $\cap$-Menger
algebras of $n$-place functions} (cf. \cite{Dudtro2}) isomorphic
to some Menger algebras of rank $n$ (cf. also \cite{DT1} and
\cite{Dudtro3}). The first abstract characterization of such
Menger $\cap$-algebras of $n$-place functions was announced
(without proof) in \cite{Trrel}. The full proof of this result
can be found in the monograph \cite{Dudtro2} (see Theorem 3.4.2).
Putting in this theorem $n=1$ we obtain the following statement.

\begin{theorem}\label{T-prpol1roda}
An algebraic system $(G,\cdot,\wedge,\sqsubset)$, where
$(G,\cdot)$ is a semigroup, $(G,\wedge)$ is a
semilattice, $\sqsubset$ is a binary relation on $G$,
is isomorphic to some projection $\cap$-semigroup of
transformations if and only if $\sqsubset$ is a left regular
and right negative quasi-order containing the semilattice order $\leq$
$($defined by $x\leq y\longleftrightarrow x\wedge y=x)$ and
the identity \eqref{dt-1} and the conditions
\begin{eqnarray}
&&\label{f-tr11} x\sqsubset x\wedge y\longrightarrow x\leq y,\\
&&\label{f-tr12} (x\wedge y) u\leq yu,\\
&&\label{f-tr13} x\sqsubset y\wedge z\mathrel{\&} x\sqsubset (z\wedge v)u
\longrightarrow x\sqsubset (y\wedge z\wedge v)u
\end{eqnarray}
are satisfied for all $x,y,z,v\in G$ and $u\in G^1$.
\end{theorem}

In the same manner, from the results proved in \cite{Dudtro2} we
obtain:

\begin{theorem}\label{T-prpol1rodob}
An algebraic system $(G,\cdot,\wedge,\sqsubset)$,
where $(G,\cdot)$ is a semigroup, $(G,\wedge)$ is a
semilattice, $\sqsubset$ is a binary relation on $G$, is isomorphic to
some projection $\cap$-semigroup of invertible transformations if
and only if it satisfies the identity
\begin{equation}\label{dt-3}
(x\wedge y)z=xz\wedge yz
\end{equation}
and all conditions of Theorem \emph{\ref{T-prpol1roda}}.
\end{theorem}

\medskip

{\bf 4.} The semiadjacency relation was first described on
semigroups of transformations. Namely, in \cite{pavl-2}
transformation semigroup of the form $(\Phi,\circ,\delta_{\Phi})$
are characterized by some essentially infinite system of elementary
axioms. Ordered and projection ordered semigroups of transformations
equipped with the semiadjacency relation are characterized in \cite{garv-trokh}.
Here we characterize $\cap$-semi\-groups of transformations of
the form $(\Phi,\circ,\cap,\delta_{\Phi})$ and
$(\Phi,\circ,\cap,\chi_{\Phi},\delta_{\Phi})$.

\begin{theorem}\label{T-trt1}
An algebraic system $(G,\cdot,\wedge,\sqsubset,\vdash)$,
where $(G,\cdot)$ is a semigroup, $(G,\wedge)$ is a
semilattice, $\sqsubset$, $\vdash$ are binary relations on $G$,
is isomorphic to some projection $\cap$-semigroup
$(\Phi,\circ,\cap,\chi_{\Phi},\delta_{\Phi})$ of transformations
if and only if all conditions of Theorem \emph{\ref{T-prpol1roda}}
are satisfied and the equivalence
\begin{equation}\label{f-tr14}
  x\vdash y\longleftrightarrow x\sqsubset xy,
\end{equation}
holds for all $x,y\in G$.
\end{theorem}

\begin{proof}
The necessity of these conditions follows from our Theorem
\ref{T-prpol1roda} and from condition (\ref{f-tr3}) of Proposition
\ref{P-trp1}.

To prove the sufficiency of these conditions we assume that
$(G,\cdot,\wedge,\sqsubset,\vdash)$ sa\-tis\-fies all conditions
mentioned in Theorem~\ref{T-trt1}. Since the conditions of
Theorem~\ref{T-prpol1roda} are satisfied,
$(G,\cdot,\wedge,\sqsubset)$ is isomorphic to a projection
$\cap$-semigroup $(\Phi,\circ,\cap,\chi_{\Phi})$ of
transformations. Let $P:G\rightarrow\Phi$ be this
isomorphism. Thus, we have
\begin{eqnarray*}
&& P(g_1g_2)=P(g_2)\circ P(g_1),\\
&& P(g_1\wedge g_2)=P(g_1)\cap P(g_2),\\
&& g_1\sqsubset g_2\longleftrightarrow (P(g_1),P(g_2))\in\hi
\end{eqnarray*}
for any $g_1,g_2\in G$. The last, together with \eqref{f-tr14},
shows that $x\vdash y$ is equivalent to $(P(x),P(xy))\in\hi$,
i.e., to $(P(x),P(y)\circ P(x))\in\hi$. Hence, by Proposition~\ref{P-trp1},
we get $(P(x),P(y))\in\delta_{\Phi}$.
So,
\[
x\vdash y\longleftrightarrow (P(x),P(y))\in\delta_{\Phi}.
\]
Therefore $P$ is an isomorphism of
$(G,\cdot,\wedge,\sqsubset,\vdash)$ onto a projection
$\cap$-semigroup $(\Phi,\circ,\cap,\hi,\delta_{\Phi})$ of
transformations.
\end{proof}

Similarly, we can prove the following theorem:

\begin{theorem}\label{T-trt1*}
An algebraic system
$(G,\cdot,\wedge,\sqsubset,\vdash)$, where $(G,\cdot)$ is a
semigroup, $(G,\wedge)$ is a semilattice, $\sqsubset$, $\vdash$
are binary relations on $G$, is isomorphic to some projection
$\cap$-semigroup $(\Phi,\circ,\cap,\hi,\delta_{\Phi})$ of
invertible transformations if and only if \eqref{f-tr14} and all
conditions of Theorem \emph{\ref{T-prpol1rodob}} are satisfied.
\end{theorem}

\medskip

{\bf 5.} Let $(G,\cdot,\wedge,\vdash)$ be an algebraic system
such that $(G,\cdot)$ is a semigroup, $(G,\wedge)$ is a
semilattice, $\vdash$ is a left ideal relation on a semigroup
$(G,\cdot)$, and the conditions \eqref{dt-1}, \eqref{f-tr11},
\eqref{f-tr12} are satisfied. For brevity, the formula $x\vdash
y\mathrel{\&}  xy\leq z$ will be denoted by $x\boxdot y\leq z$.
By $G^{*}$ we denote the set $G\cup\{e\}$, where $e$ is some element,
not belonging to $G$, such that (by  definition) $ee=e$, $e\leq e$,
$e\vdash e$, $ex = xe = x$ and $x\vdash e$ for all $x\in G$.

\begin{definition}\label{D-d2}
\rm A subset $H\subset G$ is called \emph{$f$-closed} if
\begin{equation}\label{f-18}
   (u\wedge v\wedge w)x\boxdot y\leq zt\,\mathrel{\&} \,(u\wedge v\in H)\mathrel{\&}
   ((v\wedge w)x\in H)\longrightarrow z\in H
\end{equation}
for any $x,y,t\in G^*$ and $z,u,v,w\in G$.
\end{definition}

The smallest $f$-closed subset of $G$ containing $H\subset G$ is
denoted by $f(H)$. It can be characterized by sets of the form
$F_n(H)$ defined as follows: $F_0(H)=H$,
\[
F_1(H) =\{z\,|\,(\exists u,v,w,x,y,t)\, \big((u\wedge
v\wedge w)x \boxdot y\leq zt\mathrel{\&} \big(u\wedge v\in
H\big)\mathrel{\&} \big((v\wedge w) x\in H\big)\big\},
\]
where $z,u,v,w\in G$, $x,y,t\in G^*,$ and $F_n(H)=F_1(F_{n-1}(H))$
for every $n>1$.

\begin{lemma}\label{lem}
$H\subset F_1(H)$ for every subset $H$ of $G$ and $H=F_1(H)$ for
every $f$-closed subset $H$ of $\,G$.
\end{lemma}
\begin{proof}
Since $G\subset G^*,$ for any $z\in H\subset G$ we have
$(z\wedge z\wedge z)e\boxdot e\leq ze$, $z\wedge
z\in H$ and $(z\wedge z)e\in H$. Hence $z\in F_1(H)$. So,
$H\subset F_1(H)$.

If $H$ is $f$-closed and $z\in F_1(H)$, then $(u\wedge
v\wedge w)x\boxdot y\leq zt$, $u\wedge v\in H$ and
$(v\wedge w)x\in H$ for some $u,v,w\in G$, $x,y,t\in G^*$.
Hence, by \eqref{f-18}, we obtain $z\in H$. So, $F_1(H)\subset H$.
This completes the proof.
\end{proof}

\begin{corollary}\label{C1}
$H=F_0(H)\subset F_1(H)\subset F_2(H)\subset F_3(H)\subset\ldots $
\ for every subset $H$ of $G$.
\end{corollary}

\begin{proposition}\label{P-p4}
\ $f(H)=\bigcup\limits_{n\in\mathbb{N}}\!\!F_n(H)$\ for any
$H\subset G$, where $\mathbb{N}$ is the set of non-negative
integers.
\end{proposition}

\noindent\emph{Proof}. Let $(u\wedge v\wedge w)x\boxdot
y\leq zt$, \ $u\wedge
v\in\!\!\bigcup\limits_{n\in\mathbb{N}}\!\!F_n(H)$ and
$(v\wedge
w)x\in\!\!\bigcup\limits_{n\in\mathbb{N}}\!\!F_n(H)$. Then
$u\wedge v\in F_{n_1}(H)$ and $(v\wedge w)x\in
F_{n_2}(H)$ for some $n_1$, $ n_2$. Let $m=\max\{n_1,n_2\}$. Since
$F_{n_i}(H)\subset F_{m}(H)$, $i=1,2$, we have
$$
(u\wedge v \wedge w) x\boxdot y\leq zt\mathrel{\&}
(u\wedge v\in F_m(H))\mathrel{\&} \big((v\wedge w) x\in
F_{m}(H)\big),
$$
which implies $z\in
F_m(H)\subset\bigcup\limits_{n\in\mathbb{N}}\!\! F_n(H)$. This
proves that $\bigcup\limits_{n\in\mathbb{N}}F_n(H)$ is $f$-closed.
Since $H\subset\bigcup\limits_{n\in\mathbb{N}}\!\! F_n(H)$
by Corollary \ref{C1}, it follows that $f(H)\subset\bigcup\limits_{n\in\mathbb{N}}\!\! F_n(H)$
by the minimality of $f(H)$. But, as it is not difficult to see,
$F_1(A)\subset F_1(B)$ for $A\subset B\subset G$, which together
with Lemma \ref{lem} implies $F_1(H)\subset F_1(f(H))= f(H)$.
Similarly, $F_2{}(H)\subset f(H)$ and so on. Hence,
$\bigcup\limits_{n\in\mathbb{N}}F_{n}(H)\subset f(H)$, which
completes the proof.
 \hfill $\Box{}$

\medskip

Using the method of mathematical induction we can prove the
following proposition.

\begin{proposition}\label{P-p5} For each subset $H$ of the
algebraic system $(G,\cdot,\wedge,\vdash)$ and each natural
number $n>1$ an element $z\in G$ belongs to $F_{n}(H)$ if and only
if the system of conditions
\begin{equation}\label{f-20}
 \hspace*{-2ex} \left.\begin{array}{c}
      (u_1\wedge v_1\wedge w_1) x_1\boxdot y_1\leq zt_1,\\
      \bigwedge\limits_{i=1}^{2^{n-1}-1}\left(\begin{array}{l}
        (u_{2i}\wedge v_{2i}\wedge w_{2i})x_{2i}\boxdot y_{2i}\leq
        (u_i\wedge v_i)t_{2i},\\
        (u_{2i+1}\wedge v_{2i+1}\wedge w_{2i+1})x_{2i+1}\boxdot y_{2i+1}
        \leq(v_i\wedge w_i)x_it_{2i+1}
     \end{array}\right),\\
     \bigwedge\limits_{i=2^{n-1}}^{2^n-1}(u_i\wedge v_i,(v_i\wedge w_i)x_i\in H)
\end{array}\right\}
\end{equation}
is satisfied for some $x_i,y_i,t_i\in G^*$ and $u_i,v_i,w_i\in G$.
\end{proposition}

\begin{proposition}\label{P-p6}
In a $\cap$-semigroup $(\Phi,\circ,\cap,\delta_{\Phi})$
of transformations for any subset $H_{\Phi}\subset\Phi$ and any
$\varphi\in f(H_{\Phi})$ we have $\bigcap\limits_{\varphi_i\in
H_{\Phi}}\!\!\!\mathrm{pr}_1\varphi_i\subset\mathrm{pr}_1\varphi$.
\end{proposition}

\begin{proof} First of all, we prove that the implication
\begin{equation}\label{f-21}
  \varphi\in F_{n}(H_{\Phi})\longrightarrow\bigcap\limits_{\varphi_i\in H_{\Phi}}
\!\!  \mathrm{pr}_1\varphi_i\subset\mathrm{pr}_1\varphi
\end{equation}
is true for every $n\geq 0$. For simplicity we denote the
intersection $\bigcap\limits_{\varphi_i\in
H_{\Phi}}\!\!\!\pr\varphi_i$ by $\mathfrak{A}$.

For $n=0$ this implication is obvious. Assume that it is true for
some $n$. We prove that it is true for $n+1$. For this, select an
arbitrary $\varphi\in F_{n+1}(H_{\Phi})$. Then, according to the
definition of $F_{n+1}(H)$, for some $h,f,g,p,s,r\in\Phi$ we have
$\,(h\circ (f\cap g\cap r),p)\in\delta_{\Phi}$,\ $p\circ h\circ
(f\cap g\cap r)\subset s\circ\varphi$,\ $f\cap g\in
F_{n}(H_{\Phi})$ and $h\circ (g\cap r)\in F_{n}(H_{\Phi})$. From
$f\cap g\in F_{n}(H_{\Phi})$ it follows
$\mathfrak{A}\subset\pr(f\cap g)$. Similarly, $h\circ (g\cap r)\in
F_{n}(H_{\Phi})$ implies $\mathfrak{A}\subset\pr(h\circ (g\cap
r))$. The condition $(h\circ (f\cap g\cap r), p)\in\delta_{\Phi}$
means that $\prr (h\circ (f\cap g\cap r))\subset\pr p$, whence,
according to \eqref{f-tr3}, we obtain $\pr (h\circ (f\cap g\cap
r))\subset\pr(p\circ h\circ (f\cap g\cap r))$. Analogously $p\circ
h\circ (f\cap g\cap r)\subset s\circ\varphi$ gives $\pr (p\circ
h\circ (f\cap g\cap r))\subset\pr (s\circ\varphi)$. So, $\pr
(h\circ (f\cap g\cap
r))\subset\pr(s\circ\varphi)\subset\pr\varphi$.

Next, from $\mathfrak{A}\subset\pr (f\cap g)$ and
$\mathfrak{A}\subset\pr (h\circ (g\cap r))$ we obtain
\begin{multline*}
\mathfrak{A}=\Delta_{\spr (f\cap g)}(\mathfrak{A})\subset\Delta_{\spr (f\cap g)}
(\pr (h\circ (g\cap r)))=\\
\pr \left((h\circ (g\cap r))\circ\Delta_{\spr (f\cap g)}\right)=\pr(h\circ (f\cap
g\cap r)),
\end{multline*}
where $\Delta_{\spr (f\cap g)}$ is the identity transformation on
the set $\pr (f\cap g)$.\footnote{This means that $\Delta_{\spr (f\cap
g)}=\{(a,a)\,|\,a\in\pr (f\cap g)\}$.}

 Thus, $\mathfrak{A}\subset\pr\varphi$. So, the implication
\eqref{f-21} is true for $n+1$. Consequently, it is true for all
$n\geq 0$.

Let now $\varphi\in f(H_{\Phi})$, then, according to Proposition
\ref{P-p4}, there is $n$ such that $\varphi\in F_{n}(H_{\Phi})$,
whence, applying \eqref{f-21}, we deduce
$\mathfrak{A}\subset\pr\varphi$.
\end{proof}

\begin{proposition}\label{P-p7}
A subset $H$ of the algebra $(G,\cdot,\wedge,\vdash)$ is
$f$-closed if and only if it satisfies following four conditions:
\begin{eqnarray}
& &\label{f-22} xy\in H\longrightarrow x\in H,\\
& &\label{f-23} g_1\vdash g_2\mathrel{\&}  g_1\in H\longrightarrow g_1g_2\in H,\\
& &\label{f-24} g_1\wedge g_2=g_1\in H\longrightarrow g_2\in H,\\
& &\label{f-25} g_1\wedge g_2\in H\mathrel{\&}  (g_2\wedge g_3)
h\in H\longrightarrow (g_1 \wedge g_2\wedge g_3) h\in H,
\end{eqnarray}
for all $x,y,g_1,g_2,g_3\in G$, $h\in G^*$.
\end{proposition}

\begin{proof}
Let $H$ be an $f$-closed subset of $G$. Then,
according to the definition,
\begin{equation}\label{f-26}
\left.\begin{array}{c} (u\wedge v\wedge w) x\vdash
y\mathrel{\&}  (u\wedge v\wedge w) xy\leq zt,\\[2pt]
 u\wedge v\in H\mathrel{\&}  (v\wedge w) x\in H
 \end{array}\right\}\longrightarrow z\in H,
\end{equation}
where $x,y,t\in G^*$, $z,u,v\in G$.

Putting in this formula $u=v=w=xy$, $x=y=e$, $t=y$ and $z=x$, we
obtain the implication
\[\left.\begin{array}{c}
(xy\wedge xy \wedge xy)e\vdash e\mathrel{\&}  (xy\wedge
xy\wedge xy) ee\leq xy,\\[2pt]
xy\wedge xy\in H\mathrel{\&}  (xy\wedge xy)e\in H
\end{array}\right\}\longrightarrow x\in H,
\]
i.e., the implication $xy\vdash e\mathrel{\&}  xy\leq xy\mathrel{\&}  xy\in
H\longrightarrow x\in H$. Since $xy\vdash e$ and $xy\leq xy$ are
always true, the above means that $xy\in H\longrightarrow x\in H$.
So, \eqref{f-22} is valid for any $f$-closed subset $H$.

In a similar way putting in \eqref{f-26} $u=v=w=g_1$, $x=e$,
$y=g_2$, $t=e$, $z=g_1g_2$,  we obtain \eqref{f-23}, putting
$u=v=w=g_1$, $x=y=t=e$, $z=g_2$, we obtain \eqref{f-24}. Next,
putting $u=g_1$, $v=g_2$, $w=g_3$, $y=e$, $z=(g_1\wedge
g_2\wedge g_3)h$, $t=e$, we get \eqref{f-25}. So, every
$f$-closed subset $H$ satisfies \eqref{f-22}--\eqref{f-25}.

Conversely, let \eqref{f-22}--\eqref{f-25} and the premise of
\eqref{f-26} be satisfied. Then, from $u\wedge v\in H$ and
$(v\wedge w)x\in H$, by \eqref{f-25}, we obtain
$(u\wedge v\wedge w)x\in H$. Next, from $(u\wedge
v\wedge w)x\vdash y$ and $(u\wedge v\wedge w)x\in
H$, applying \eqref{f-23}, we deduce $(u\wedge v\wedge
w)xy\in H$. Now, from $(u\wedge v\wedge w)xy\leq zt$ and
$(u\wedge v\wedge w)xy\in H$, according to \eqref{f-24},
we conclude $zt\in H$, whence, by \eqref{f-22}, we get $z\in H$.
So, \eqref{f-22}--\eqref{f-25} imply \eqref{f-26}. Hence $H$ is
$f$-closed.
\end{proof}

Consider on $(G,\cdot,\wedge,\vdash)$ the binary relation
$\chi_0$ defined as follows:
\[
(g_1,g_2)\in\chi_{_0}\longleftrightarrow \,g_2\in\langle
g_1\rangle,
\]
where $\langle g_1\rangle$ is the $f$-closure $f(\{g_1\})$ of $\{g_1\}$.

\begin{proposition}\label{P-p8} The relation $\chi_{_0}$ is
the smallest relation among left regular and right negative
quasi-orders $\sqsubset$ defined on $(G,\cdot,\wedge,\delta)$,
containing the natural order $\leq$ on the semilattice $(G,\wedge)$,
and satisfying the conditions
\begin{eqnarray}
& &\label{f-27} g_1\vdash g_2\longrightarrow g_1\sqsubset
g_1g_2,\\
& &\label{f-28} g\sqsubset g_1\wedge g_2\mathrel{\&}  g\sqsubset
(g_2\wedge g_3)x\longrightarrow g\sqsubset (g_1\wedge
g_2\wedge g_3)x,
\end{eqnarray}
where $g_1,g_2,g_3,g\in G$, $x\in G^*$.
\end{proposition}

\begin{proof} First  observe that $\chi_{_0}$ is a quasi-order.
Since $x\in\langle x\rangle$, we have $(x,x)\in\chi_{_0}$ for any
$x\in G$. If $(x,y)\in\chi_{0}$ and $(y,z)\in\chi_{_0}$, then
$y\in\langle x\rangle$ and $z\in\langle y\rangle$, whence we
obtain $z\in\langle x\rangle$. So, $(x,z)\in\chi_{_0}$. Thus
$\chi_{_0}$ is reflexive and transitive, i.e., it is a quasi-order
on $G$.

Now, we show that
\begin{equation}\label{f-29a}
   zf(H)\subset f(zH)
\end{equation}
for any subset $H\subset G$ and any element $z\in G$.

In view of Proposition \ref{P-p4} it is sufficient to prove that
\begin{equation}\label{f-30a}
   zF_{n}(H)\subset F_{n}(zH)
\end{equation}
for each $n=0,1,2,\ldots$

For $n=0$ this inclusion is obvious. Suppose that it is true for
some $n>0$. To prove it for $n+1$ consider an arbitrary element
$g\in zF_{n+1}(H)$. Then $g=zg_1$ for some $g_1\in F_{n+1}(H)$.
But for any $g_1\in F_{n+1}(H)$ there are $u,v,w, x,y,t$ such that
\[
(u\wedge v\wedge w)x\vdash y ,\ (u\wedge
v\wedge w) xy\leq g_1t, \ u\wedge v\in F_{n}(H) \ {\rm
and } \ (v\wedge w) x\in F_{n}(H).
\]
Since the relation $\vdash$ is left ideal, from $(u\wedge
v\wedge w) x\vdash y$ we get $z(u\wedge v\wedge w)
x\vdash y$, whence, using \eqref{dt-1}, we obtain $(zu\wedge
zv\wedge zw)x\vdash y$. Next, from $(u\wedge
v\wedge w)xy\leq g_1t$, in view of the stability of $\leq$,
we deduce $z(u\wedge v\wedge w) xy\leq zg_1t$, which, by
\eqref{dt-1}, gives $(zu\wedge zv\wedge zw)xy\leq
zg_1t$.

Clearly, $u\wedge v\in F_{n}(H)$ and $(v\wedge w)x\in
F_{n}(H)$ imply $z(u\wedge v)\in zF_{n}(H)$ and
$z(v\wedge w)x\in zF_{n}(H)$, whence, by \eqref{dt-1}, we
obtain $zu\wedge zv\in zF_{n}(H)$ and $(zv\wedge zw)x\in
zF_{n}(H)$. This, according to \eqref{f-30a} and our assumption on
$n$, gives $zu\wedge zv,\ (zv\wedge zw)x\in F_{n}(zH)$.

So, we have $(zu\wedge zv\wedge zw)x\vdash y$,\
$(zu\wedge zv\wedge zw)xy\leq zg_1t$,\ $zu\wedge
zv\in F_{n}(zH)$ and $(zv\wedge zw)x\in F_{n}(zH)$, which
implies $zg_1\in F_{n+1}(zH)$, i.e., $g\in F_{n+1}(zH)$. This
means that the implication \eqref{f-30a} is true for $n+1$.
Consequently, it is true for every $n\geq 0$.

Therefore
\[
zf(H)=z\Big(\bigcup\limits_{n\in\mathbb{N}}F_{n}(H)\Big)\subset\bigcup\limits_{n\in\mathbb{N}}
F_{n}(zH)=f(zH),
\]
which proves \eqref{f-29a}.

Now, we can verify the properties of the relation $\chi_{_0}$. At
first observe that if $(x,y)\in\chi_{_0}$, then $y\in\langle
x\rangle$, whence, by \eqref{f-29a}, we obtain $zy\in z\langle
x\rangle\subset\langle zx\rangle$. Thus, $zy\in\langle zx\rangle$,
i.e., $(zx,zy)\in\chi_{_0}$. So, $\chi_{_0}$ is left regular.
Moreover, if $(z,xy)\in\chi_{_0}$, i.e., $xy\in\langle z\rangle$,
then $x\in\langle z\rangle$, by \eqref{f-22}. Hence
$(z,x)\in\chi_{_0}$. So, $\chi_{_0}$ is right negative. It also
contains the natural order of a semilattice $(G,\wedge)$.
Indeed, if $g_1\leq g_2$, then
\[
(g_1\wedge g_1\wedge g_1) e\vdash e, \ \ (g_1
\wedge g_1\wedge g_1) ee\leq g_2e \ \ {\rm and } \ \ g_1
\wedge g_1=(g_1\wedge g_1) e=g_1,
\]
whence it follows $g_2\in F_1(\{g_1\})$. As $F_1(\{g_1\})\subset
f(\{g_1\})=\langle g_1\rangle$, by Proposition~\ref{P-p4}, from
the above it follows $g_2\in\langle g_1\rangle$, i.e.,
$(g_1,g_2)\in\chi_{_0}$.

To prove \eqref{f-27} assume that $g_1\vdash g_2$ for some
$g_1,g_2\in G$. Then
\[
(g_1\wedge g_1\wedge g_1)e\vdash g_2\;\mathrel{\&}  \
(g_1\wedge g_1\wedge g_1) eg_2 \leq g_1g_2e \ \mathrel{\&}  \
g_1\wedge g_1=(g_1\wedge g_1)e=g_1,
\]
whence $g_1g_2\in F_1(\{g_1\})\subset\langle g_1\rangle$. Thus
$(g_1,g_1g_2)\in\chi_{_0}$. So, $\chi_{_0}$ satisfies
\eqref{f-27}.

Now, let $(g,g_1\wedge g_2)\in\chi_{_0}$ and $(g,(g_2
\wedge g_3)x)\in\chi_{_0}$, i.e., $g_1\wedge
g_2\in\langle g\rangle$ and $(g_2\wedge g_3)x\in\langle
g\rangle$. Since $\langle g\rangle$ is $f$-closed, the last,
according to \eqref{f-25}, implies $(g_1\wedge g_2\wedge
g_3)x\in\langle g\rangle$, whence $(g,(g_1\wedge
g_2\wedge g_3)x)\in\chi_{_0}$. This means that $\chi_{_0}$
satisfies \eqref{f-28}. So, $\chi_{_0}$ satisfies all conditions
mentioned in Proposition \ref{P-p8}.

Observe that for every $g\in G$ and an arbitrary relation $\sqsubset$
satisfying all conditions of Proposition \ref{P-p8}, the
$\sqsubset$-filter $\chi\langle g\rangle$ containing the element $g$
is an $f$-closed set\footnote{We denote by $\chi\langle g\rangle$ the
set of all elements $h\in G$ such that $g\sqsubset h$.}. Indeed,
if the premise of \eqref{f-18} is satisfied, then for some
$x,y,t\in G^*$ and $z,u,v,w\in G$ we have $(u\wedge
v\wedge w)x\vdash y$,\ $(u\wedge v\wedge w)xy\leq
zt$,\ $u\wedge v\in\chi\langle g\rangle$ and $(v\wedge
w) x\in\chi\langle g\rangle$. From $(u\wedge v\wedge
w)x\vdash y$, by \eqref{f-27}, it follows $(u\wedge
v\wedge w)x\sqsubset (u\wedge v\wedge w)xy$. From
$(u\wedge v\wedge w)xy\leq zt$, we obtain $(u\wedge v\wedge
w)xy\sqsubset zt$ since the semilattice order $\leq$ is contained
in $\sqsubset$. Thus, $(u\wedge v\wedge
w)x\sqsubset(u\wedge v\wedge w)xy\sqsubset zt$ which, by
the transitivity of $\sqsubset$, gives $(u\wedge v\wedge w)
x\sqsubset zt$. Since $\sqsubset$ is right negative, the last implies
$(u\wedge v\wedge w)x\sqsubset z$. But by \eqref{f-28},
from $u\wedge v,\ (v\wedge w)x\in\chi\langle g\rangle$
it follows $(u\wedge v\wedge w) x\in\chi\langle
g\rangle$. Thus $g\sqsubset (u\wedge v\wedge w)x$.
Consequently, $g\sqsubset z$, i.e., $z\in\chi\langle g\rangle$.
So, $\chi\langle g\rangle$ satisfies the condition \eqref{f-18}.
Hence, it is $f$-closed.

Now let $(g_1,g_2)\in\chi_{_0}$, i.e., $g_2\in\langle g_1\rangle$.
As $\chi\langle g_1\rangle$ is $f$-closed and $g_1\in\chi\langle
g_1\rangle$, we have
\[
\langle g_1\rangle =f(\{g_1\})\subset f(\chi\langle g_1\rangle)
=\chi\langle g_1 \rangle,
\]
which means that $g_2\in\chi\langle g_1\rangle$, i.e.,
$g_1\sqsubset g_2$. So, $\chi_{_0}$ is contained in $\sqsubset$.
Thus $\chi_{_0}$ is the smallest relation satisfying all conditions
of Proposition~\ref{P-p8}.
\end{proof}

Using the above results we can prove the main theorem of this
paper.

\begin{theorem}\label{T-t2} An algebraic system
$(G,\cdot,\wedge,\vdash)$, where $(G,\cdot)$ is a semigroup,
$(G,\wedge)$ is a semilattice, $\vdash$ is a binary relation
on $G$, is isomorphic to a $\cap$-semi\-group
$(\Phi,\circ,\cap,\delta_{\Phi})$ of transformations if and only
if $\delta$ is a left ideal relation on a semigroup $(G,\cdot)$,
and the following conditions
\begin{eqnarray}
& &\label{f-29} x(y\wedge z)=xy\wedge xz,\\
& &\label{f-30} (x\wedge y) z\leq yz,\\
& &\label{f-31} (x,x\wedge y)\in\chi_{_0}\longrightarrow x
\leq y,\\
& &\label{f-32} (x,xy)\in\chi_{_0}\longrightarrow x\vdash y
\end{eqnarray}
are satisfied for all $x,y,z\in G$, where $\chi_{_0}$ is as in
Proposition \emph{\ref{P-p8}} and $\leq$ is the natural order of the
semilattice $(G,\wedge)$.
\end{theorem}
\begin{proof} {\sc Necessity}. Let $P:G\to\Phi$ be an isomorphism between $(G,\cdot,\wedge,\vdash)$
and $(\Phi,\circ,\cap,\delta_{\Phi})$. Then, obviously,
$P(xy)=P(y)\circ P(x)$, $P(x\wedge y)=P(x)\cap P(y)$ and
$x\vdash y\longleftrightarrow (P(x),P(y))\in\delta_{\Phi}$
for all $x,y\in G$. Consider the relation $\chi_P$ defined on $G$
as follows:
\[
(g_1,g_2)\in\chi_P\longleftrightarrow (P(g_1),P(g_2))\in\hi,
\]
i.e., $(g_1, g_2)\in\chi_P\longleftrightarrow\pr P(g_1)\subset\pr
P(g_2)$. Clearly, the algebraic system
$(G,\cdot,\wedge,\chi_P,\vdash)$ is isomorphic to
$(\Phi,\circ,\cap,\hi,\delta_{\Phi})$, hence it satisfies all the
conditions of Theorem \ref{T-trt1}. Therefore, according to
Proposition \ref{P-p8}, $\chi_{_0}\subset\chi_P$.

Conditions \eqref{f-29}, \eqref{f-30} follow from Theorem~\ref{T-trt1}
and Proposition~\ref{P-trp3}. To prove \eqref{f-31}
let $(x,x\wedge y)\in\chi_P$. Then $(P(x),P(x)\cap
P(y))\in\chi_{\Phi}$, i.e., $\pr P(x)\subset\pr (P(x)\cap
P(y))\subset\pr P(x)$, whence $\pr P(x)=\pr (P(x)\cap P(y))$. But
$P(x)\cap P(y)\subset P(x)$, so, the last implies $P(x)\cap
P(y)=P(x)$, which is equivalent to $P(x\wedge y)=P(x)$.
Consequently, $x\wedge y=x$, i.e.\ $x\leq y$. This proves
\eqref{f-31}.

Now let $(x,xy)\in\chi_{_0}$, then $(x,xy)\in\chi_{P}$. Therefore
$(P(x),P(y)\circ P(x))\in\chi_{\Phi}$, i.e., $\pr P(x)\subset
\pr(P(y)\circ P(x))$, whence, according to \eqref{f-tr3}, we
obtain $(P(x),P(y))\in\delta_{\Phi}$. Thus, $x\vdash y$.
Condition \eqref{f-32} is proved.

{\sc Sufficiency}. Let $(G,\cdot,\wedge,\vdash)$ satisfy all
the conditions of our theorem. Proposition \ref{P-p8} shows that
the relation $\chi_{_0}$ is a left regular and right negative
quasi-order satisfying \eqref{f-27} and \eqref{f-28} and
containing the natural order of the semilattice
$(G,\wedge)$. From \eqref{f-27} and \eqref{f-32} it follows
that $x\vdash y\longleftrightarrow (x,xy)\in\chi_{_0}$ for
all $x,y\in G$. Thus, $(G,\cdot,\wedge,\chi_{_0},\vdash)$
satisfies all conditions of Theorem~\ref{T-trt1} whence it is
isomorphic to some $\cap$-semigroup
$(\Phi,\circ,\cap,\hi,\delta_{\Phi})$ of transformations. Hence
$(G,\cdot,\wedge,\vdash)$ is isomorphic to
$(\Phi,\circ,\cap,\delta_{\Phi})$.
\end{proof}

The implication \eqref{f-31} can be rewritten to the form:
\[
   x\wedge y\in\langle x\rangle\longrightarrow x\leq y,
\]
which, in view of Proposition \ref{P-p4}, is equivalent to
\[
x\wedge y\in\bigcup\limits_{n\in\mathbb{N}}F_{n}(\{x\})
\longrightarrow x\leq y,
\]
i.e., to the system of conditions $(A'_n)_{n\in\mathbb{N}}$, where
\[
A'_n\colon \; x\wedge y\in F_{n}(\{x\})\longrightarrow x \leq
y.
\]
According to Proposition \ref{P-p5} the condition $A'_n$ is
equivalent to
\[
A_n\colon\;\mathfrak{M}_n(x\wedge y,\{x\})\longrightarrow x\leq y,
\]
where $\mathfrak{M}_n(x\wedge y,\{x\})$ means the formula
\eqref{f-20} with $z=x\wedge y$ and $H=\{x\}$. So, the
implication \eqref{f-31} is equivalent to the system of conditions
$(A_n)_{n\in\mathbb{N}}$.

Similarly, we can show that the implication \eqref{f-32} is
equivalent to the system of conditions $(B_n)_{n\in\mathbb{N}}$,
where
\[
B_n\colon\,\mathfrak{M}_n(xy,\{x\})\longrightarrow x\vdash y.
\]

Thus we obtain the following

\begin{theorem}\label{T-t3}
An algebraic system $(G,\cdot,\wedge,\vdash)$, where
$(G,\cdot)$ is a semigroup, $(G,\wedge)$ is a semilattice,
$\vdash$ is a binary relations on $G$, is isomorphic to a
$\cap$-semigroup $(\Phi,\circ,\cap,\delta_{\Phi})$ of
transformations if and only if $\vdash$ is a left ideal relation
on a semigroup $(G,\cdot)$, and the system of conditions
$(A_n)_{n\in\mathbb{N}}$, $(B_n)_{n\in\mathbb{N}}$, as well as the
conditions \eqref{f-29} and \eqref{f-30} are satisfied.
\end{theorem}

\medskip

{\bf 6.} From Theorem \ref{T-trt1*} and results proved in Section~5
it follows that addding \eqref{dt-3} to the conditions of
Theorem \ref{T-trt1}, \ref{T-t2} and \ref{T-t3} we obtain an
abstract characterization of $\cap$-semigroups
$(\Phi,\circ,\cap,\delta_{\Phi})$ of invertible transformations,
as well as projection $\cap$-semigroups
$(\Phi,\circ,\cap,\hi,\delta_{\Phi})$ of invertible
transformations.

\begin{minipage} {60mm}
\begin{flushleft}
Dudek~W. A.\\ Institute of Mathematics and Computer Science\\
Wroclaw University of Technology\\ 50-370 Wroclaw\\ Poland\\
Email: dudek@im.pwr.wroc.pl
\end{flushleft}
\end{minipage}
\hfill
\begin{minipage} {60mm}
\begin{flushleft}
 Trokhimenko~V. S.\\ Department of Mathematics\\
 Pedagogical University\\
 21100 Vinnitsa\\
 Ukraine\\
 Email: vtrokhim@sovamua.com
\end{flushleft}
\end{minipage}

\end{document}